\documentclass{article}
\usepackage{latexsym,amsfonts,amssymb,amsthm,amsmath}

\newtheorem{theorem}{Theorem}

\newtheorem{proposition}[theorem]{Proposition}

\input{tcilatex}
\begin{document}

\title{Geometric View of Measurement Errors}
\author{DIARMUID O'DRISCOLL AND DONALD E. RAMIREZ}
\date{}
\maketitle

\begin{abstract}
\noindent The slope of the best fit line from minimizing the sum of the
squared oblique errors is the root of a polynomial of degree four. This
geometric view of measurement errors is used to give insight into the
performance of various slope estimators for the measurement error model
including an adjusted fourth moment estimator introduced by Gillard and Iles
(2005) to remove the jump discontinuity in the estimator of Copas (1972).
The polynomial of degree four is associated with a minimun deviation
estimator. A simulation study compares these estimators showing improvement
in bias and mean squared error.
\end{abstract}

\noindent \textsc{Department of Mathematics and Computer Studies, Mary
Immaculate College, Limerick, Ireland} \textit{E-mail address:
diarmuid.odriscoll@mic.ul.ie} \newline
\newline
\textsc{Department of Mathematics, University of Virginia, Charlottesville,
VA 22904} \textit{E-mail address: der@virginia.edu} \newline
\newline

\textbf{Keywords }Oblique errors; Measurement errors; Maximum likelihood
estimation; Moment estimation.\bigskip

\textbf{Mathematics Subject Classification} 62J05; 62G05.\bigskip

\section{Introduction}

With ordinary least squares OLS($y|x$) regression, we have data $%
\{(x_{1},Y_{1}|X=x_{1}),...,(x_{n},Y_{n}|X=x_{n})\}$ and we minimize the sum
of the squared \textit{vertical errors} to find the \emph{best-fit} line $%
y=h(x)=\beta _{0}+\beta _{1}x.$ With OLS($y|x$) it is assumed that the
independent or causal variable is measured without error. The measurement
error model has wide interest with many applications. See for example
Carroll \textit{et al.} (2006) and Fuller (1987). The comparison of
measurements by two analytical methods in clinical chemistry is often based
on regression analysis. There is no causal or independent variable in this
type of analysis. The most frequently used method to determine any
systematic difference between two analytical methods is OLS($y|x$) which has
several shortcomings when both measurement sets are subject to error.
Linnet(1993) states that \textquotedblleft it is rare that one of the
(measurement) methods is without error.\textquotedblright\ Linnet(1999)
further states that \textquotedblleft\ A systematic difference between two
(measurement) methods is identified if the estimated intercept differs
significantly from zero (constant difference) or if the slope deviates
significantly from 1 (proportional difference).\textquotedblright\ Our paper
concentrates on how to determine whether or not there is a significant
difference between two measurement instruments using a Monte Carlo
simulation; that is, we concentrate our studies about a true slope of 1. As
in the regression procedure of Deming (1943), to account for both sets of
errors, we determine a fit so that both the squared vertical and the squared
horizontal errors will be minimized. The squared vertical errors are the
squared distances from $(x,y)$ to $(x,h(x))$ and the squared horizontal
errors are the squared distances from $(x,y)$ to $(h^{-1}(y),y)$. As a
compromise, we will consider oblique errors. All of the estimated regression
models we consider (including the geometric mean and perpendicular methods)
are contained in the parametrization (with $0\leq \lambda \leq 1$) of the
line from $(x,h(x))$ to $(h^{-1}(y),y)$.

We review the Oblique Error Method in Section 2. In Section 3, we review the
geometric mean and perpendicular error models. In Section 4, we show how the
geometric mean slope is a natural estimator for the slope in the measurement
error (error-in-variables) model. Section 5 shows a relationship between the
maximum likelihood estimator in the measurement error model and the
geometric mean estimator. We give a case study to illustrate the effects
that erroneous assumptions for the ratio of variance of errors can have on
the maximum likelihood estimators. Section 6 discusses a fourth moment
estimator and shows a circular relationship to the maximum likelihood
estimator. Section 7 develops a minimum deviation estimator derived by
minimizing Equation (\ref{standard}) in Section 2 with respect to $\lambda $
for fixed $\beta _{1}$. Section 8 contains our Monte Carlo simulations where
we compare these estimators. Supporting Maple worksheets are available from
the link http://people.virginia.edu/\symbol{126}der/ODriscoll\_Ramirez/.

\section{Minimizing Squared Oblique Errors}

>From the data point $(x_{i},y_{i})$ to the fitted line $y=h(x)=\beta
_{0}+\beta _{1}x$ the vertical length is $a_{i}=\left\vert y_{i}-\beta
_{0}-\beta _{1}x_{i}\right\vert, $ the horizontal length is $%
b_{i}=\allowbreak \left\vert x_{i}-(y_{i}-\beta _{0})/\beta _{1}\right\vert
=\allowbreak \left\vert (\beta _{1}x_{i}-y_{i}+\beta _{0})/\beta
_{1}\right\vert =\allowbreak \left\vert a_{i}/\beta _{1}\right\vert $ and
the perpendicular length is $h_{i}=a_{i}/\sqrt{1+\beta _{i}^{2}}.$ With
standard notation, $S_{xx}=\sum_{i=1}^{n}(x_{i}-\overline{x})^{2},$ $%
S_{yy}=\sum_{i=1}^{n}(y_{i}-\overline{y})^{2},$ $S_{xy}=\sum_{i=1}^{n}(x_{i}-%
\overline{x})(y_{i}-\overline{y})$ with the correlation $\rho =S_{xy}/\sqrt{%
S_{xx}S_{yy}}$.

For the oblique length from $(x_{i},y_{i})$ to $(h^{-1}(y_{i})+\lambda
(x_{i}-h^{-1}(y_{i})),y_{i}+\lambda (h(x_{i})-y_{i})),$ the horizontal error
is $(1-\lambda )b_{i}=(1-\lambda )a_{i}/\left\vert\beta _{1}\right\vert$ and
the vertical error is $\lambda a_{i}.$ The sum of squared horizontal,
respectively vertical, errors are given by $SSE_{h}(\beta _{0},\beta
_{1},\lambda)=\left( \sum_{i=1}^{n}a_{i}^{2}\right) /\beta _{1}^{2}$ and $%
SSE_{v}(\beta_{0},\beta _{1},\lambda )=\sum_{i=1}^{n}a_{i}^{2}$. In a
comprehensive paper by Riggs \textit{et al.} (1978), the authors place great
emphasis on the importance of equations being dimensionally correct, since
it is from these equations that the slope estimators are derived. In
particular the authors state that: ``It is a poor method indeed whose
results depend upon the particular units chosen for measuring the variables
... and that invariance under linear transformations is equivalent to
requiring the method to be dimensionally correct.'' So that our equation is
dimensionally correct we consider

\begin{equation}
SSE_{o}(\beta _{0},\beta _{1},\lambda )=(1-\lambda )^{2}\frac{SSE_{h}}{%
\widetilde{\sigma }_{\delta }^{2}}+\lambda ^{2}\frac{SSE_{v}}{\widetilde{%
\sigma }_{\tau }^{2}}  \label{nor}
\end{equation}%
where $\{\widetilde{\sigma }_{\delta }^{2},\widetilde{\sigma }_{\tau }^{2}\}$
are Madansky's moment estimators of the variance in the horizontal,
respectively vertical, directions. In Section 4, we show that this is
equivalent to using 
\begin{equation}
SSE_{o}(\beta _{0},\beta _{1},\lambda )=(1-\lambda
)^{2}S_{yy}SSE_{h}+\lambda ^{2}S_{xx}SSE_{v}.  \label{standard}
\end{equation}%
Similar to that shown in O'Driscoll, Ramirez and Schmitz (2008), the
solution of $\partial SSE_{o}/\partial \beta _{0}=0$ is given by $\beta _{0}=%
\overline{y}-\beta _{1}\overline{x}$ and hence 
\begin{equation*}
SSE_{o}(\beta _{0},\beta _{1},\lambda )=\left( (1-\lambda )^{2}S_{yy}/\beta
_{1}^{2}+\lambda ^{2}S_{xx}\right) (S_{yy}-2\beta _{1}S_{xy}+\beta
_{1}^{2}S_{xx}).
\end{equation*}%
The solutions of $\partial SSE_{o}/\partial \beta _{1}=0$ are then the roots
of the fourth degree polynomial equation in $\beta _{1}$, namely%
\begin{equation}
P_{4}(\beta _{1})=\lambda ^{2}\sqrt{\dfrac{S_{xx}}{S_{yy}}}\dfrac{S_{xx}}{%
S_{yy}}\beta _{1}^{4}-\lambda ^{2}\dfrac{S_{xx}}{S_{yy}}\rho \beta
_{1}^{3}+(1-\lambda )^{2}\rho \beta _{1}-(1-\lambda )^{2}\sqrt{\dfrac{S_{yy}%
}{S_{xx}}}=0.  \label{P4}
\end{equation}%
With $\lambda =1$ we recover the minimum squared vertical errors with
estimated slope $\beta _{1}^{ver}$, and with $\lambda =0$ we recover the
minimum squared horizontal errors with estimated slope $\beta _{1}^{hor}$.

For each fixed $\lambda \in \lbrack 0.1]$, there corresponds $\beta _{1}\in
\lbrack \beta _{1}^{ver},\beta _{1}^{hor}]$ which satisfies Equation (\ref%
{P4}), and conversely, for each fixed $\beta _{1}\in \lbrack \beta
_{1}^{ver},\beta _{1}^{hor}]$, there corresponds $\lambda \in \lbrack 0.1]$
such that minimizing the sum of the squared oblique errors has estimated
slope $\beta _{1}$. We measure the angle $\theta _{\lambda }$ of the oblique
projection associated with $\lambda $ using the line segments $(x,y)$ to $%
(x,h(x))$ and $(x,h(x))$ to $(h^{-1}(y),y)$. When the slope $\beta _{1}$ is
close to one, for $\lambda $ near one we anticipate $\theta _{\lambda }$ to
be near 45$^{\circ }$ and for $\lambda $ is close to zero we anticipate $%
\theta _{\lambda }$ to be near 135$^{\circ }.$ The angles are computed from
the Law of Cosines.

A similar argument to that of O'Driscoll \textit{et al.} (2008) shows that $%
P_{4}(\beta _{1})$ has exactly two real roots, one positive and one negative
with the global minimum being the positive (respectively negative) root
corresponding to the sign of $S_{xy}$. Riggs \textit{et al.} (1978) in
Equation (119) also noted the role of the roots of a similar quartic
equation in determining the slope estimators.

\section{Minimizing Squared Perpendicular and Squared Geometric Mean Errors}

The perpendicular error model dates back to Adcock (1878) who introduced it
as a procedure for fitting a straight line model to data with error measured
in both the $x$ and $y$ directions. For squared perpendicular errors Adcock
minimized $SSE_{per}(\beta _{0},\beta _{1})=\sum_{i=1}^{n}a_{i}^{2}/(1+\beta
_{1}^{2})$ with solutions $\beta _{0}^{per}=\overline{y}-\beta _{1}^{per}%
\overline{x}$ and 
\begin{equation}
\beta _{1}^{per}=\dfrac{(S_{yy}-S_{xx})\pm \sqrt{%
(S_{yy}-S_{xx})^{2}+4S_{xy}^{2}}}{2S_{xy}},  \label{b1p}
\end{equation}%
provided $S_{xy}\neq 0$. However, in this case, the equation which minimizes 
$SSE_{per}(\beta _{0},\beta _{1})$ is dimensionally incorrect unless $x$ and 
$y$ are measured in the same units.

For squared geometric mean errors, we minimize $SSE_{gm}(\beta _{0},\beta
_{1})\allowbreak =\sum_{i=1}^{n}\left( \sqrt{\lvert a_{i}b_{i}\rvert }%
\right) ^{2}\allowbreak =\sum_{i=1}^{n}a_{i}^{2}/\lvert \beta _{1}\rvert $
with solutions $\beta _{0}^{gm}=\overline{y}-\beta _{1}^{gm}\overline{x}$
and 
\begin{equation}
\beta _{1}^{gm}=\pm \sqrt{S_{yy}/S_{xx}}.  \label{b1gm}
\end{equation}

\begin{proposition}
The geometric mean estimator has oblique parameter $\lambda =1/2.$
\end{proposition}

Proof: For $\beta _{1}^{gm}=\sqrt{S_{yy}/S_{xx}}$,we solve the quadratic
equation $P_{4}(\beta _{1}^{gm})=0$ for $\lambda $. This equation reduces to
a linear equation whose root is $\lambda =1/2.\square $

\section{Measurement Error Model and Second Moment Estimation}

We now consider the measurement error (errors-in-variables) model as
follows. In this paper it is assumed that $X$ and $Y$ are random variables
with respective finite variances $\sigma _{X}^{2}$ and $\sigma _{Y}^{2}$,
finite fourth moments and have the linear functional relationship $Y=\beta
_{0}+\beta _{1}X$. The observed data $\{(x_{i},y_{i}),1\leq i\leq n\}$ are
subject to error by $x_{i}=X_{i}+\delta _{i}$ and $y_{i}=Y_{i}+\tau _{i}$
where it is also assumed that $\delta $ is $N(0,\sigma _{\delta }^{2})$ and $%
\tau $ is $N(0,\sigma _{\tau }^{2})$. In our simulation studies we will use
an exponential distribution for $X$.

It is well known, in a measurement error model, that the expected value for $%
\beta _{1}^{ver}$ is attenuated to zero by the attenuating factor $\sigma
_{X}^{2}/(\sigma _{\delta }^{2}+\sigma _{X}^{2})$, called the reliability
ratio by Fuller (1987). Similarly the expected value for $\beta _{1}^{hor}$
is amplified to infinity by the amplifying factor $(\sigma _{Y}^{2}+\sigma
_{\tau }^{2})/\sigma _{Y}^{2}$. Thus for the measurement error model, when
both the vertical and horizontal models are reasonable, a compromise
estimator such as the geometric mean estimator $\beta _{1}^{gm}$ is hoped to
have improved efficiency. However, Lindley and El-Sayyad (1968) proved that
the expected value of $\beta _{1}^{gm}$ is biased unless $\sigma _{\tau
}^{2}/\sigma _{Y}^{2}=\sigma _{\delta }^{2}/\sigma _{X}^{2}$.

Madansky's moment estimators for $\{\sigma _{\delta }^{2},\sigma _{\tau
}^{2}\}$ are%
\begin{eqnarray}
\widetilde{\sigma }_{\delta }^{2} &=&\frac{S_{xx}}{n}-\frac{S_{xy}}{n\beta
_{1}},  \label{sigmas} \\
\widetilde{\sigma }_{\tau }^{2} &=&\frac{S_{yy}}{n}-\frac{\beta _{1}S_{xy}}{n%
}.  \notag
\end{eqnarray}%
If $\sigma _{\delta }$ is known or can be approximated, Madansky used the
first of the equations in Equation (\ref{sigmas}) to derive an estimator for 
$\beta _{1}$ but Riggs \textit{et al.} (1978) in Figure 5 produced an
example for which this estimator performs poorly. In general $\sigma
_{\delta }$ is not known and concentration focuses on estimating the true
error ratio $\sigma _{\tau }^{2}/\sigma _{\delta }^{2}$. Our Monte Carlo
simulation illustrates how poor estimates for the error ratio may lead to
large biases using the MLE estimator. It would be interesting to determine
if there is a slope estimator $\beta _{1}$ that is a fixed point for $\beta
_{1}=\widetilde{\sigma }_{\tau }(\beta _{1})/\widetilde{\sigma }_{\delta
}(\beta _{1})$. This can be achieved with the geometric mean estimator $%
\beta _{1}^{gm}$.

\begin{proposition}
$\beta _{1}^{gm}$ is a fixed point of the ratio function $\beta _{1}=\dfrac{%
\widetilde{\sigma }_{\tau }(\beta _{1})}{\widetilde{\sigma }_{\delta }(\beta
_{1})}.$
\end{proposition}

Proof: Rewrite $n\widetilde{\sigma }_{\delta }^{2}=S_{xx}-S_{xy}\beta
_{1}^{gm}=S_{xx}-S_{xy}\sqrt{S_{xx}/S_{yy}}$ and $n\widetilde{\sigma }_{\tau
}^{2}=S_{yy}-S_{xy}\beta _{1}^{gm}=S_{yy}-S_{xy}\sqrt{S_{yy}/S_{xx}}$, from
which $\widetilde{\sigma }_{\tau }/\widetilde{\sigma }_{\delta }=\sqrt{%
S_{yy}/S_{xx}}$.$\square $

We return to the assertion made in Section 1. A natural standardized weighed
average for the oblique model is shown in Equation (\ref{nor}) and using the
fixed point solution of Proposition 2 in this equation yields the equivalent
model given in Equation (\ref{standard}).

\section{The Maximum Likelihood Estimator}

If the ratio of the error variances $\kappa =\sigma _{\tau }^{2}/\sigma
_{\delta }^{2}$ is assumed finite, then Madansky (1959), among others,
showed that the maximum likelihood estimator for the slope is 
\begin{equation}
\beta _{1}^{mle}=\frac{(S_{yy}-\kappa S_{xx})+\sqrt{(S_{yy}-\kappa
S_{xx})^{2}+4\kappa \rho ^{2}S_{xx}S_{yy}}}{2\rho \sqrt{S_{xx}S_{yy}}}
\label{b1mle}
\end{equation}

It also follows that if $\kappa =1$ in Equation (\ref{b1mle}) then the MLE
(often called the Deming Regression estimator) is equivalent to the
perpendicular estimator, $\beta _{1}^{per}$. Conversely, if the MLE is $%
\beta _{1}^{per}$ then $\kappa =1$. In the particular case where $%
S_{xx}=S_{yy}$ then $\beta _{1}^{per}$ has a $\lambda $ value of $0.5$. We
note that $S_{yy}/S_{xx}$ is a good estimator of $\sigma _{y}^{2}/\ \sigma
_{x}^{2}$, but in general, it is not a good estimator of the error ratio $%
\kappa =\sigma _{\tau }^{2}/\ \sigma _{\delta }^{2}$. In Section 6, we
discuss a moment estimator $\widetilde{\kappa }$ for $\kappa .$

The MLE is a function of $\{\rho ,\kappa ,S_{xx}/S_{yy}\}$. Table 1 gives
the corresponding $\beta _{1}^{mle}$ value for typical values. For fixed $%
\{\kappa ,\rho \}$, the values for $\beta _{1}^{mle}$ in each column of
Table 1 decrease. As expected, with $\kappa =1$ and $S_{xx}=S_{yy}$, the
maximum likelihood estimator agrees with the geometric mean estimator, both
being equal to $1.00$.\bigskip

\begin{center}
\textbf{Table 1}

Values for $\beta _{1}^{mle}$ for typical $\{\rho ,\kappa ,S_{xx}/S_{yy}\}$

\begin{tabular}{lllllllllllll}
\hline
$\kappa =$ & $0.500$ & $0.500$ & $0.500$ & $0.500$ & $1.000$ & $1.000$ & $%
1.000$ & $1.000$ & $2.000$ & $2.000$ & $2.000$ & $2.000$ \\ 
$\rho =$ & $0.200$ & $0.400$ & $0.600$ & $0.800$ & $0.200$ & $0.400$ & $0.600
$ & $0.800$ & $0.200$ & $0.400$ & $0.600$ & $0.800$ \\ \hline
$S_{xx}/S_{yy}=1/2$ & $5.396$ & $2.828$ & $2.016$ & $1.632$ & $3.799$ & $%
2.219$ & $1.750$ & $1.535$ & $1.414$ & $1.414$ & $1.414$ & $1.414$ \\ 
$S_{xx}/S_{yy}=1$ & $2.686$ & $1.569$ & $1.237$ & $1.086$ & $1.000$ & $1.000$
& $1.000$ & $1.000$ & $0.372$ & $0.638$ & $0.808$ & $0.921$ \\ 
$S_{xx}/S_{yy}=2$ & $0.707$ & $0.707$ & $0.707$ & $0.707$ & $0.263$ & $0.451$
& $0.571$ & $0.651$ & $0.185$ & $0.354$ & $0.496$ & $0.613$ \\ \hline
\end{tabular}%
\bigskip 
\end{center}

In Table 2 we record the corresponding obliqueness parameter $\lambda $ for
the maximum likelihood model for these typical values. Small values near $0$
support OLS($x|y$), denoted by $\beta _{1}^{hor}$, and large values near $1$
support OLS($y|x$), denoted by $\beta _{1}^{ver}$. For fixed $\{\kappa ,\rho
\}$, the values for the obliqueness parameter $\lambda $ in each column of
Table 2 increase indicating the model moves from $\beta _{1}^{hor}$ towards $%
\beta _{1}^{ver}$. With $\kappa =S_{yy}/S_{xx}$, $\beta _{1}^{mle}=\beta
_{1}^{gm}$ as shown by the cells of Table 2 with $\lambda =0.500$.\bigskip

\begin{center}
\textbf{Table 2}

Values for $\lambda $ for typical $\{\rho ,\kappa ,S_{xx}/S_{yy}\}$

\begin{tabular}{cllllllllllll}
\hline
$\kappa =$ & $0.500$ & $0.500$ & $0.500$ & $0.500$ & $1.000$ & $1.000$ & $%
1.000$ & $1.000$ & $2.000$ & $2.000$ & $2.000$ & $2.000$ \\ 
$\rho =$ & $0.200$ & $0.400$ & $0.600$ & $0.800$ & $0.200$ & $0.400$ & $0.600
$ & $0.800$ & $0.200$ & $0.400$ & $0.600$ & $0.800$ \\ \hline
$S_{xx}/S_{yy}=1/2$ & $0.033$ & $0.111$ & $0.197$ & $0.273$ & $0.089$ & $%
0.223$ & $0.316$ & $0.375$ & $0.500$ & $0.500$ & $0.500$ & $0.500$ \\ 
$S_{xx}/S_{yy}=1$ & $0.089$ & $0.223$ & $0.316$ & $0.375$ & $0.500$ & $0.500$
& $0.500$ & $0.500$ & $0.911$ & $0.777$ & $0.684$ & $0.625$ \\ 
$S_{xx}/S_{yy}=2$ & $0.500$ & $0.500$ & $0.500$ & $0.500$ & $0.911$ & $0.776$
& $0.684$ & $0.625$ & $0.967$ & $0.889$ & $0.803$ & $0.727$ \\ \hline
\end{tabular}%
\textit{\bigskip }
\end{center}

The Madansky's moment estimators $\{\widetilde{\sigma }_{\delta }^{2},%
\widetilde{\sigma }_{\tau }^{2}\}$ depend on the choice of $\beta _{1}$. In
Table 3, we record the effect of varying slopes on the moments and their
ratio when computable.\bigskip

\begin{center}
\textbf{Table 3}

Error ratios for Madansky's moment estimators for varying $\beta _{1}$

\begin{tabular}{lccc}
\hline
& $\widetilde{\sigma }_{\delta }^{2}$ & $\widetilde{\sigma }_{\tau }^{2}$ & $%
\frac{\widetilde{\sigma }_{\tau }^{2}}{\widetilde{\sigma }_{\delta }^{2}}$
\\ \hline
$\beta _{1}^{ver}$ & $0$ & $\frac{1-\rho ^{2}}{n}S_{yy}$ & $\infty $ \\ 
$\beta _{1}^{hor}$ & $\frac{1-\rho ^{2}}{n}S_{xx}$ & $0$ & $0$ \\ 
$\beta _{1}^{gm}$ & $\frac{1-\rho }{n}S_{xx}$ & $\frac{1-\rho }{n}S_{yy}$ & $%
\frac{S_{yy}}{S_{xx}}$ \\ 
$\beta _{1}^{per}$ & $\frac{1}{2}\frac{S_{xx}+S_{yy}-\sqrt{%
(S_{xx}-S_{yy})^{2}+4\rho ^{2}S_{xx}S_{yy}}}{n}$ & $\frac{1}{2}\frac{%
S_{xx}+S_{yy}-\sqrt{(S_{xx}-S_{yy})^{2}+4\rho ^{2}S_{xx}S_{yy}}}{n}$ & $1$
\\ 
$\beta _{1}^{mle}$ & $\frac{1}{2}\frac{S_{xx}+\frac{S_{yy}}{\kappa }-\sqrt{%
(S_{xx}-\frac{S_{yy}}{\kappa })^{2}+4\rho ^{2}S_{xx}\frac{S_{yy}}{\kappa }}}{%
n}$ & $\frac{1}{2}\frac{\kappa S_{xx}+S_{yy}-\sqrt{(\kappa
S_{xx}-S_{yy})^{2}+4\rho ^{2}\kappa S_{xx}S_{yy}}}{n}$ & $\kappa $ \\ \hline
\end{tabular}%
\bigskip 
\end{center}

In the next section, we introduce a second moment estimator for $\kappa $
and a fourth moment estimator for $\beta _{1}.$

\section{Fourth Moment Estimation}

When $\kappa $ is unknown, Solari (1969) showed that the maximum likelihood
estimator for the slope $\beta _{1}$ does not exist, as the maximum
likelihood surface has a saddle point at the critical value. Earlier Lindley
and El-Sayyad (1968) suggested, in this case, that the maximum likelihood
method fails as the estimator would be the geometric mean estimator which
converges to the wrong value. Sprent (1970) pointed out the result of Solari
does not imply that the maximum likelihood principle has failed, but rather
that the likelihood surface has no maximum value at the critical value.

Copas (1972) offered some advice for using the maximum likelihood method. He
assumed the data has rounding-off errors in the observations which allows
for an approximated likelihood function to be used, and that this
approximated likelihood function is bounded. His estimator for the slope has
the rule%
\begin{equation*}
\beta _{1}^{cop}=\left\{ 
\begin{array}{c}
\beta _{1}^{ver}\text{ \ \ if }\sum y_{i}^{2}<\sum x_{i}^{2} \\ 
\beta _{1}^{hor}\text{ \ \ if }\sum y_{i}^{2}>\sum x_{i}^{2}%
\end{array}%
\right. ,
\end{equation*}%
so the ordinary least squares estimators are used depending on the whether $%
|\beta _{1}^{gm}|~<1$ or $|\beta _{1}^{gm}|~>1.$

The Copas estimator is \textit{not} continuous in the data as a small change
in data can switch the direction of the inequality $\sum y_{i}^{2}<\sum
x_{i}^{2}$ which will cause a jump discontinuity in the estimator $\beta
_{1}^{cop}.$ To achieve continuity in the data, we adjust the range of the
fourth moment estimator $\beta _{1}^{mom}$ described in Gillard and Iles
(2005) to account for admissible values for $\{\sigma _{\delta }^{2},\sigma
_{\tau }^{2}\}$. See also Gillard and Iles (2010).

The basic second moment estimators for $\widetilde{\sigma }_{\delta }^{2}$
and $\widetilde{\sigma }_{\tau }^{2}$ are shown in Equation (\ref{sigmas}).
Since variances must be positive, we have the admissible range for the
moment estimator for $\widetilde{\beta }_{1}$ as%
\begin{equation}
\beta _{1}^{ver}=\frac{S_{xy}}{S_{xx}}<\widetilde{\beta }_{1}<\frac{S_{yy}}{%
S_{xy}}=\beta _{1}^{hor}.  \label{beta mom range}
\end{equation}%
Set $S_{xxxy}=\sum (x_{i}-\overline{x})^{3}(y_{i}-\overline{y})$ and
similarly for $S_{xyyy}.$ Following Gillard and Iles (2005), from Equations
(22) and (24)), the fourth moment equations of interest are%
\begin{eqnarray}
\frac{S_{xxxy}}{n} &=&\widetilde{\beta }\widetilde{\mu }_{4}+3\widetilde{%
\beta }\widetilde{\sigma }^{2}\widetilde{\sigma }_{\delta }^{2}
\label{22 & 24} \\
\frac{S_{xyyy}}{n} &=&\widetilde{\beta }^{3}\widetilde{\mu }_{4}+3\widetilde{%
\beta }\widetilde{\sigma }^{2}\widetilde{\sigma }_{\varepsilon }^{2}  \notag
\end{eqnarray}%
with $\widetilde{\mu }_{4}$ denoting the fourth central moment for the
underlying distribution of $X.$ The four equations from (\ref{sigmas}) and (%
\ref{22 & 24}) allow for a moment solution for $\beta _{1}$ as 
\begin{equation}
\widetilde{\beta }_{1}=\sqrt{\frac{S_{xyyy}-3S_{xy}S_{yy}}{%
S_{xxxy}-3S_{xy}S_{xx}}.}  \label{beta Gillard}
\end{equation}%
In our simulation study, $\widetilde{\beta }_{1}$ was well-defined around $%
99\%$ of the time. If the radicand is negative, we recommend using the
geometric mean estimator.

To satisfy Equation (\ref{beta mom range}) we define $\beta _{1}^{mom}$ as%
\begin{equation}
\beta _{1}^{mom}=\left\{ 
\begin{array}{l}
\beta _{1}^{ver}\text{ \ \ if }\widetilde{\beta }_{1}\leq \beta _{1}^{ver}
\\ 
\widetilde{\beta }_{1}\text{ \ \ \ \ if }\beta _{1}^{ver}\leq \widetilde{%
\beta }_{1}\leq \beta _{1}^{hor} \\ 
\beta _{1}^{hor}\text{ \ \ if }\widetilde{\beta }_{1}\geq \beta _{1}^{hor}%
\end{array}%
\right. .  \label{beta mom}
\end{equation}%
This is a Copas-type estimator with the moment estimator $\widetilde{\beta }%
_{1}$ used to "smooth" out the jump discontinuity inherent in the Copas
estimator. We next study the circular relationship between this moment
estimator and the maximum likelihood estimator with fixed $\kappa $.

We will define the moment estimator $\kappa (\beta _{1})$ as a function of $%
\beta _{1}$, then use this value to compute $\beta _{1}^{mle}(\kappa )$ as a
function of $\kappa .$ Finally, we note that $\beta _{1}^{mle}(\kappa (%
\widetilde{\beta }_{1}))=\widetilde{\beta }_{1},$ showing the circular
relationship between the estimators $\{\widetilde{\beta }_{1},\beta
_{1}^{mle}\}.$ Thus our moment estimator also has the functional form of the
maximum likelihood estimator with fixed $\kappa .$

Set $\widetilde{\kappa }(\widetilde{\beta }_{1})=\widetilde{\sigma }_{\tau
}^{2}/\widetilde{\sigma }_{\delta }^{2}$ so%
\begin{equation}
\widetilde{\kappa }(\widetilde{\beta }_{1})=\frac{S_{yy}-\widetilde{\beta }%
_{1}\rho \sqrt{S_{xx}S_{yy}}}{S_{xx}-\rho /\widetilde{\beta }_{1}\sqrt{%
S_{xx}S_{yy}}}.  \label{kappa tilda}
\end{equation}%
We use $\widetilde{\kappa }(\widetilde{\beta }_{1})$ in Equation (\ref{b1mle}%
) to determine $\beta _{1}^{mle}(\widetilde{\kappa }(\widetilde{\beta }%
_{1})).$ As $\widetilde{\beta }_{1}\rightarrow \beta _{1}^{hor}$ the
numerator in Equation (\ref{kappa tilda}) tends to zero so $\widetilde{%
\kappa }(\widetilde{\beta }_{1})\rightarrow 0$ and $\beta _{1}^{mle}(%
\widetilde{\kappa }(\widetilde{\beta }_{1}))\rightarrow \beta _{1}^{hor}$;
similarly as $\widetilde{\beta }_{1}\rightarrow \beta _{1}^{ver}$ the
denominator in Equation (\ref{kappa tilda}) tends to zero so $\widetilde{%
\kappa }(\widetilde{\beta }_{1})\rightarrow \infty $ and $\beta _{1}^{mle}(%
\widetilde{\kappa }(\widetilde{\beta }_{1}))\rightarrow \beta _{1}^{ver}.$ A
stronger result is given in the following Proposition.

\begin{proposition}
For each $\beta _{1}$, $\beta _{1}^{mle}(\widetilde{\kappa }(\beta
_{1}))=\beta _{1}$ and in particular $\beta _{1}^{mle}(\widetilde{\kappa }(%
\widetilde{\beta }_{1}))=\widetilde{\beta }_{1}.$
\end{proposition}

Proof: In Equation (\ref{b1mle}) solve $\beta _{1}^{mle}(\kappa )=\beta _{1}$
for $\kappa =\kappa _{0}$, and then check that $\kappa _{0}$ is the same as
in Equation (\ref{kappa tilda}).$\square $

An example helps to demonstrate the smoothing achieved with the moment
estimator $\beta _{1}^{mom}.$ Assume $\{\rho
=0.5,S_{xx}=1,S_{xxxy}=10,S_{xyyy}=5\}.$ Equation (\ref{beta mom range})
requires that $0.13029\leq S_{yy}\leq 1.31862.$ As $S_{yy}$ varies over the
admissible values for $S_{yy}$, $\widetilde{\kappa }(\widetilde{\beta }_{1})$
varies over $[0,\infty ]\,$and $\widetilde{\beta }_{1}$ varies over $[\beta
_{1}^{ver},\beta _{1}^{hor}]$ and $\beta _{1}^{mle}(\widetilde{\kappa }(%
\widetilde{\beta }_{1}))=\widetilde{\beta }_{1}$, a surprising result$.$%
\bigskip 

\begin{center}
\textbf{Table 4}

Slope Estimates with $\{\rho =0.5,S_{xx}=1,S_{xxxy}=10,S_{xyyy}=5\}$ 
\begin{equation*}
\begin{tabular}{cccccc}
\hline
$S_{yy}$ & $\beta _{1}^{ver}$ & $\widetilde{\beta }_{1}$ & $\beta _{1}^{hor}$
& $\widetilde{\kappa }(\widetilde{\beta }_{1})$ & $\beta _{1}^{mle}$ \\ 
\hline
$0.1303$ & $0.1805$ & $0.7219$ & $0.7219$ & $0.0000$ & $0.7219$ \\ 
$0.2000$ & $0.2236$ & $0.7222$ & $0.8944$ & $0.0558$ & $0.7222$ \\ 
$0.4000$ & $0.3164$ & $0.7145$ & $1.2649$ & $0.3123$ & $0.7145$ \\ 
$0.6000$ & $0.3873$ & $0.6977$ & $1.5492$ & $0.7412$ & $0.6977$ \\ 
$0.8000$ & $0.4472$ & $0.6734$ & $1.7889$ & $1.4850$ & $0.6734$ \\ 
$1.0000$ & $0.5000$ & $0.6417$ & $2.0000$ & $3.0760$ & $0.6417$ \\ 
$1.2000$ & $0.5477$ & $0.6020$ & $2.1909$ & $9.6582$ & $0.6020$ \\ 
$1.3186$ & $0.5742$ & $0.5742$ & $2.2966$ & $\infty $ & $0.5741$ \\ \hline
\end{tabular}%
\end{equation*}%
\bigskip 
\end{center}

\section{Minimum Deviation Estimation}

>From Section 1 with fixed $\beta _{1}$ the solution of $\partial
SSE_{o}/\partial \lambda =0$ is given by $\lambda =S_{yy}/(S_{yy}+\beta
_{1}^{2}S_{xx})$. Substituting $\beta _{1}^{mom}$ for $\beta _{1}$ in this
result for $\lambda $ produces a Minimum Deviation type estimator which we
denote by $\beta _{1}^{md}$, with $\beta _{1}^{ver}\leq \beta _{1}^{md}\leq
\beta _{1}^{hor}$. Our simulation studies will support the efficiency of
this Minimum Deviation estimator.

Riggs \textit{et al.} (1978) observes for $S_{yy}/S_{xx}=1$ that the
geometric mean estimator is unbiased when $\kappa =1$, negatively biased
when $\kappa <1,$ and positively biased when $\kappa >1$. The authors also
state that \textquotedblleft no one method of estimating $\beta _{1}$ is the
best method under all circumstances.\textquotedblright\ To determine the
efficiency of these estimators we conduct a Monte Carlo simulation in the
next section.

\section{Monte Carlo Simulation}

Our Monte Carlo simulation uses $X$ with an exponential distribution with
mean $\mu _{X}=10$ (and $\sigma _{X}=10)$ and $Y=X$ so $\beta _{1}=1$ and $%
\beta _{0}=0$. Both $X$ and $Y$ are subject to errors $\sigma _{\delta }^{2}$%
, respectively $\sigma _{\tau }^{2}$ where $(\sigma _{\delta }^{2}$, $\sigma
_{\tau }^{2})$ $\in $ $\{1,4,9\}$$\times $$\{1,4,9\}$. The sample size $n$
is chosen as $100$.

The first simulation, with the number of replications $R=$ $100$, summarized
in Table 5, reports on the bias in the MLE estimator in using a misspecified
value of $\kappa $. For $(\sigma _{\delta }^{2}$, $\sigma _{\tau }^{2})$ $%
\in $ $\{1,4,9\}$$\times $$\{1,4,9\},$ $\kappa $ ranges with ratios from $%
1:1 $ to $9:1$ The true error ratios of $\kappa $ are recorded in the first
row and the assumed error ratios $\kappa ^{\#}$ which are used to compute $%
\beta _{1}^{mle}$ are recorded in the first column, both in ascending
order.\bigskip

\begin{center}
\textbf{Table 5}

Percentage Bias of MLE estimator for the assumed ratios $\kappa ^{\#}$ for
varying values of $\kappa =\sigma _{\tau }^{2}/\sigma _{\delta }^{2}$

$\{\beta _{1}=1,\beta _{0}=0,n=100,R=100\}$

\begin{tabular}{cccccccccc}
\hline
$\{\kappa ^{\#},\kappa \}$ & $1:9$ & $1:4$ & $4:9$ & $1:1$ & $4:4$ & $9:9$ & 
$9:4$ & $4:1$ & $9:1$ \\ \hline
$1:9$ & \multicolumn{1}{r}{$0.166$} & \multicolumn{1}{r}{$0.502$} & 
\multicolumn{1}{r}{$2.164$} & \multicolumn{1}{r}{$0.870$} & 
\multicolumn{1}{r}{$3.663$} & \multicolumn{1}{r}{$7.995$} & 
\multicolumn{1}{r}{$8.723$} & \multicolumn{1}{r}{$3.592$} & 
\multicolumn{1}{r}{$9.282$} \\ 
$1:4$ & \multicolumn{1}{r}{$-0.914$} & \multicolumn{1}{r}{$-0.012$} & 
\multicolumn{1}{r}{$0.811$} & \multicolumn{1}{r}{$0.666$} & 
\multicolumn{1}{r}{$2.807$} & \multicolumn{1}{r}{$6.087$} & 
\multicolumn{1}{r}{$7.351$} & \multicolumn{1}{r}{$3.067$} & 
\multicolumn{1}{r}{$8.265$} \\ 
$4:9$ & \multicolumn{1}{r}{$-2.066$} & \multicolumn{1}{r}{$-0.564$} & 
\multicolumn{1}{r}{$-0.643$} & \multicolumn{1}{r}{$0.445$} & 
\multicolumn{1}{r}{$1.878$} & \multicolumn{1}{r}{$3.999$} & 
\multicolumn{1}{r}{$5.838$} & \multicolumn{1}{r}{$2.496$} & 
\multicolumn{1}{r}{$7.137$} \\ 
$1:1$ & \multicolumn{1}{r}{$-4.067$} & \multicolumn{1}{r}{$-1.541$} & 
\multicolumn{1}{r}{$-3.184$} & \multicolumn{1}{r}{$0.051$} & 
\multicolumn{1}{r}{$0.218$} & \multicolumn{1}{r}{$0.266$} & 
\multicolumn{1}{r}{$3.083$} & \multicolumn{1}{r}{$1.467$} & 
\multicolumn{1}{r}{$5.058$} \\ 
$4:4$ & \multicolumn{1}{r}{$-4.067$} & \multicolumn{1}{r}{$-1.541$} & 
\multicolumn{1}{r}{$-3.184$} & \multicolumn{1}{r}{$0.051$} & 
\multicolumn{1}{r}{$0.218$} & \multicolumn{1}{r}{$0.266$} & 
\multicolumn{1}{r}{$3.083$} & \multicolumn{1}{r}{$1.467$} & 
\multicolumn{1}{r}{$5.058$} \\ 
$9:9$ & \multicolumn{1}{r}{$-4.067$} & \multicolumn{1}{r}{$-1.541$} & 
\multicolumn{1}{r}{$-3.184$} & \multicolumn{1}{r}{$0.051$} & 
\multicolumn{1}{r}{$0.218$} & \multicolumn{1}{r}{$0.266$} & 
\multicolumn{1}{r}{$3.083$} & \multicolumn{1}{r}{$1.467$} & 
\multicolumn{1}{r}{$5.058$} \\ 
$9:4$ & \multicolumn{1}{r}{$-5.957$} & \multicolumn{1}{r}{$-2.495$} & 
\multicolumn{1}{r}{$-5.590$} & \multicolumn{1}{r}{$-0.342$} & 
\multicolumn{1}{r}{$-1.417$} & \multicolumn{1}{r}{$-3.330$} & 
\multicolumn{1}{r}{$0.338$} & \multicolumn{1}{r}{$0.437$} & 
\multicolumn{1}{r}{$2.936$} \\ 
$4:1$ & \multicolumn{1}{r}{$-6.956$} & \multicolumn{1}{r}{$-3.016$} & 
\multicolumn{1}{r}{$-6.856$} & \multicolumn{1}{r}{$-0.561$} & 
\multicolumn{1}{r}{$-2.310$} & \multicolumn{1}{r}{$-5.230$} & 
\multicolumn{1}{r}{$-1.161$} & \multicolumn{1}{r}{$-0.136$} & 
\multicolumn{1}{r}{$1.748$} \\ 
$9:1$ & \multicolumn{1}{r}{$-7.840$} & \multicolumn{1}{r}{$-3.489$} & 
\multicolumn{1}{r}{$-7.973$} & \multicolumn{1}{r}{$-0.763$} & 
\multicolumn{1}{r}{$-3.119$} & \multicolumn{1}{r}{$-6.899$} & 
\multicolumn{1}{r}{$-2.513$} & \multicolumn{1}{r}{$-0.663$} & 
\multicolumn{1}{r}{$0.657$} \\ \hline
\end{tabular}%
\bigskip 
\end{center}

As expected, the values for $\widetilde{\kappa }=\kappa $ show the smallest
bias, and in each column for a given $\kappa $ the bias shows that the
estimated slope moves from over estimating the true value to under
estimating the true value of $\beta _{1}=1.$ This was anticipated since for $%
\widetilde{\kappa }$ near zero the maximum likelihood estimator favors $%
\beta _{1}^{hor}$ which over estimates $\beta _{1;}$ and correspondingly,
for $\widetilde{\kappa }$ near one the maximum likelihood estimator favors $%
\beta _{1}^{ver}$ which under estimates $\beta .$

We conducted a second large scale Monte Carlo simulation study with $R=1000$
to demonstrate the improvement in the adjusted fourth moment estimator $%
\beta _{1}^{mom}$ over the Copas estimator which has a jump discontinuity.
Simulations for other slope estimators have been reported by Hussin (2004).
We used an exponential distribution for $X$ with $\mu _{X}=10,$ and set $%
\beta _{1}=1$ and $\beta _{1}=0.$ The values for the error standard
deviations were $(\sigma _{\delta }$, $\sigma _{\tau })$ $\in $ $\{1,2,3,4\}$%
$\times $$\{1,2,3,4\},$ the sample size was $n=100$ and the number of
replications $R=1000.$ We report in Tables 6, 7, and 8 the MSE and the Bias
for the estimators $\{\beta _{1}^{ver},\beta _{1}^{hor},\beta
_{1}^{per},\beta _{1}^{gm},\beta _{1}^{mom},\beta _{1}^{md}\}$ for $(\sigma
_{\delta }$, $\sigma _{\tau })\in \{(1,2),(1,3),(1,4)\}$. Similar results
hold for $(\sigma _{\delta }$, $\sigma _{\tau })\in \{(2,1),(3,1),(4,1)\}$
Note that in each case the adjusted fourth moment estimator $\beta _{1}^{mom}
$ is more efficient than the Copas estimator. To see this we compare the
pairs of values (MSE, Bias) in the three tables. For $\beta _{1}^{mom}$
these are $\{(1.001,-0.830),(2.786,-1.807),(5.717,-2.813)\}$ and for Copas
these are $\{(2.378,-2.410),(8.769,-7.347),(23.018,-13.848)\}.$ In practice,
the researcher may not know which of $\{\sigma _{\tau }^{2},\sigma _{\delta
}^{2}\}$ is larger. If he does, then he may choose either of $\{\beta
_{1}^{ver},\beta _{1}^{hor}\}$ with $\beta _{1}^{hor}$ favored when $\sigma
_{\delta }^{2}$ is much bigger than $\sigma _{\tau }^{2}$. A fairer
comparison is to use OLS($y|x$) and OLS($x|y$) each $50\%$ of the time. Thus
in the Tables we report the average for the MSE and the average of the
absolute deviation of the biases for the two OLS estimators. These average
(MSE, Bias) values from the tables are $%
\{(1.336,2.518),(4.847,4.831),(12.46,7.858)\}$ showing the improved
efficiency of $\beta _{1}^{mom}$. As anticipated the minimum deviation
estimator $\beta _{1}^{ml}$ achieves further improvement in reduction of
(MSE, Bias) with values $\{(0.646,-1.336),(2.309,-3.584),(5.578,-6.288)\}.$%
\bigskip 

\begin{center}
\textbf{Table 6}

$X$ is $Exp(10),\beta _{1}=1,\beta _{0}=0,R=1000,n=100$

$(\sigma _{\tau }=1,\sigma _{\delta }=2)$

$OLS^{\ast }$ reports average MSE and average absolute Bias for $\{\beta
_{1}^{ver},\beta _{1}^{hor}\}$

$%
\begin{tabular}{lcccc}
\hline
& $MSE\ast 10^{-3}$ & $\%Bias$ & $\lambda $ & $\theta _{\lambda }$ \\ \hline
$\beta _{1}^{ver}$ & \multicolumn{1}{r}{$2.001$} & \multicolumn{1}{r}{$-3.843
$} & \multicolumn{1}{r}{$1.000$} & \multicolumn{1}{r}{$46.12$} \\ 
$OLS^{\ast }$ & \multicolumn{1}{r}{$1.336$} & \multicolumn{1}{r}{$2.518$} & 
\multicolumn{1}{r}{NA} & \multicolumn{1}{r}{NA} \\ 
$\beta _{1}^{hor}$ & \multicolumn{1}{r}{$0.670$} & \multicolumn{1}{r}{$1.193$%
} & \multicolumn{1}{r}{$0.000$} & \multicolumn{1}{r}{$136.12$} \\ 
$\beta _{1}^{per}$ & \multicolumn{1}{r}{$0.688$} & \multicolumn{1}{r}{$-1.396
$} & \multicolumn{1}{r}{$0.507$} & \multicolumn{1}{r}{$89.99$} \\ 
$\beta _{1}^{gm}$ & \multicolumn{1}{r}{$0.653$} & \multicolumn{1}{r}{$-1.360$%
} & \multicolumn{1}{r}{$0.500$} & \multicolumn{1}{r}{$90.78$} \\ 
$\beta _{1}^{mom}$ & \multicolumn{1}{r}{$1.001$} & \multicolumn{1}{r}{$-0.830
$} & \multicolumn{1}{r}{$0.339$} & \multicolumn{1}{r}{$108.27$} \\ 
$\beta _{1}^{cop}$ & \multicolumn{1}{r}{$2.378$} & \multicolumn{1}{r}{$-2.410
$} & \multicolumn{1}{r}{$0.651$} & \multicolumn{1}{r}{$74.47$} \\ 
$\beta _{1}^{md}$ & \multicolumn{1}{r}{$0.646$} & \multicolumn{1}{r}{$-1.336$%
} & \multicolumn{1}{r}{$0.497$} & \multicolumn{1}{r}{$91.06$} \\ \hline
\end{tabular}%
$\textit{\bigskip }

\textbf{Table 7}

$X$ is $Exp(10),\beta _{1}=1,\beta _{0}=0,R=1000,n=100$

$(\sigma _{\tau }=1,\sigma _{\delta }=3)$

$OLS^{\ast }$ reports average MSE and average absolute Bias for $\{\beta
_{1}^{ver},\beta _{1}^{hor}\}$

$%
\begin{tabular}{lcccc}
\hline
& $MSE\ast 10^{-3}$ & $\%Bias$ & $\lambda $ & $\theta _{\lambda }$ \\ \hline
$\beta _{1}^{ver}$ & \multicolumn{1}{r}{$8.370$} & \multicolumn{1}{r}{$-8.459
$} & \multicolumn{1}{r}{$1.000$} & \multicolumn{1}{r}{$47.53$} \\ 
$OLS^{\ast }$ & \multicolumn{1}{r}{$4.847$} & \multicolumn{1}{r}{$4.831$} & 
\multicolumn{1}{r}{NA} & \multicolumn{1}{r}{NA} \\ 
$\beta _{1}^{hor}$ & \multicolumn{1}{r}{$1.324$} & \multicolumn{1}{r}{$1.203$%
} & \multicolumn{1}{r}{$0.000$} & \multicolumn{1}{r}{$137.53$} \\ 
$\beta _{1}^{per}$ & \multicolumn{1}{r}{$2.688$} & \multicolumn{1}{r}{$-3.954
$} & \multicolumn{1}{r}{$0.520$} & \multicolumn{1}{r}{$89.60$} \\ 
$\beta _{1}^{gm}$ & \multicolumn{1}{r}{$2.423$} & \multicolumn{1}{r}{$-3.760$%
} & \multicolumn{1}{r}{$0.500$} & \multicolumn{1}{r}{$92.19$} \\ 
$\beta _{1}^{mom}$ & \multicolumn{1}{r}{$2.786$} & \multicolumn{1}{r}{$-1.807
$} & \multicolumn{1}{r}{$0.318$} & \multicolumn{1}{r}{$110.94$} \\ 
$\beta _{1}^{cop}$ & \multicolumn{1}{r}{$8.769$} & \multicolumn{1}{r}{$-7.347
$} & \multicolumn{1}{r}{$0.848$} & \multicolumn{1}{r}{$58.14$} \\ 
$\beta _{1}^{md}$ & \multicolumn{1}{r}{$2.309$} & \multicolumn{1}{r}{$-3.584$%
} & \multicolumn{1}{r}{$0.490$} & \multicolumn{1}{r}{$93.196$} \\ \hline
\end{tabular}%
$\bigskip 

\textbf{Table 8}

$X$ is $Exp(10),\beta _{1}=1,\beta _{0}=0,R=1000,n=100$

$(\sigma _{\tau }=1,\sigma _{\delta }=4)$

$OLS^{\ast }$ reports average MSE and average absolute Bias for $\{\beta
_{1}^{ver},\beta _{1}^{hor}\}$

\begin{tabular}{lcccc}
\hline
& $MSE\ast 10^{-3}$ & $\%Bias$ & $\lambda $ & $\theta _{\lambda }$ \\ \hline
$\beta _{1}^{ver}$ & \multicolumn{1}{r}{$22.791$} & \multicolumn{1}{r}{$%
-14.376$} & \multicolumn{1}{r}{$1.000$} & \multicolumn{1}{r}{$49.43$} \\ 
$OLS^{\ast }$ & \multicolumn{1}{r}{$12.46$} & \multicolumn{1}{r}{$7.858$} & 
\multicolumn{1}{r}{NA} & \multicolumn{1}{r}{NA} \\ 
$\beta _{1}^{hor}$ & \multicolumn{1}{r}{$2.134$} & \multicolumn{1}{r}{$1.339$%
} & \multicolumn{1}{r}{$0.000$} & \multicolumn{1}{r}{$139.43$} \\ 
$\beta _{1}^{per}$ & \multicolumn{1}{r}{$7.406$} & \multicolumn{1}{r}{$-7.480
$} & \multicolumn{1}{r}{$0.539$} & \multicolumn{1}{r}{$89.95$} \\ 
$\beta _{1}^{gm}$ & \multicolumn{1}{r}{$6.242$} & \multicolumn{1}{r}{$-6.880$%
} & \multicolumn{1}{r}{$0.500$} & \multicolumn{1}{r}{$94.08$} \\ 
$\beta _{1}^{mom}$ & \multicolumn{1}{r}{$5.717$} & \multicolumn{1}{r}{$-2.813
$} & \multicolumn{1}{r}{$0.286$} & \multicolumn{1}{r}{$114.51$} \\ 
$\beta _{1}^{cop}$ & \multicolumn{1}{r}{$23.018$} & \multicolumn{1}{r}{$%
-13.848$} & \multicolumn{1}{r}{$0.950$} & \multicolumn{1}{r}{$52.71$} \\ 
$\beta _{1}^{md}$ & \multicolumn{1}{r}{$5.578$} & \multicolumn{1}{r}{$-6.288$%
} & \multicolumn{1}{r}{$0.480$} & \multicolumn{1}{r}{$96.04$} \\ \hline
\end{tabular}%
\bigskip 
\end{center}

\section{Summary}

We have modified the fourth moment estimator of the slope from Gillard and
Iles (2005) to show how to remove the jump discontinuity in the estimator
given by Copas (1972). We show how the moment estimators $\{\beta _{1}^{mom},%
\widetilde{\sigma }_{\delta }^{2},\widetilde{\sigma }_{\tau }^{2}\}$ can be
used to determine an MLE estimator which surprisingly is the original moment
estimator of the slope. Our simulations support our claim that both
\thinspace $\{\beta _{1}^{mom},\beta _{1}^{md}\}$ are more efficient than
the average of the OLS estimators.

\section{References}

Adcock, R. J. (1878). A problem in least-squares. The Analyst
5:53-54.\bigskip

Carroll, R. J., Ruppert, D., Stefanski, L. A., Crainiceanu, C. M. (2006). 
\textit{Measurement Error in Nonlinear Models - A Modern Perspective}, 
\textit{Second Edition}. Boca Raton: Chapman \& Hall/CRC.\bigskip

Copas, J. (1972). The likelihood surface in the linear functional
relationship problem. Journal of the Royal Statistical Society. Series B
(Methodological) 34:274-278.\bigskip

Deming, W. E. (1943). \textit{Statistical Adjustment of Data.} New York:
Wiley. \bigskip

Fuller, W.A. (1987). \textit{Measurement Error Models. }New York:\textit{\ }%
Wiley.\bigskip

Gillard J.,Iles T. (2005). Method of moments estimation in linear regression
with errors in both variables. Cardiff University School of Mathematics
Technical Report.\bigskip

Gillard, J., Iles T. (2009). Methods of fitting straight lines where both
variables are subject to measurement error. Current Clinical Pharmacology%
\textit{\ }4:164-171.\bigskip

Hussin, A.G. (2004). Numerical comparisons for various estimators of slope
parameters for unreplicated linear functional model.\bigskip Matematika%
\textit{\ }20:19-30.

Lindley, D., El-Sayyad, M. (1968). The Bayesian estimation of a linear
functional relationship. Journal of the Royal Statistical Society Series B
(Methodological) 30:190-202\bigskip

Linnet, Kristian (1993). Evaluation of regression procedures for methods
comparison studies. Clinical Chemistry\textbf{\ }39:432-432.\bigskip

Linnet, Kristian (1999). Necessary sample size for method comparison studies
based on regression analysis Clinical Chemistry\textbf{\ }45:882-894.\bigskip

Madansky, A. (1959). The fitting of straight Lines when both variables are
subject to error. Journal of American Statistical Association
54:173-205.\bigskip

O'Driscoll, D., Ramirez, D., Schmitz, R. (2008). Minimizing oblique errors
for robust estimation. Irish. Math. Soc. Bulletin\textit{\ }62\textit{:71-78}%
.\bigskip

Riggs, D., Guarnieri, J., Addelman, S. (1978). Fitting straight lines when
both variables are subject to error. Life Sciences\textbf{\ }%
22:1305-1360.\bigskip

Sprent, P. (1970) \ The saddle point of the likelihood surface for a linear
functional relationship. Journal of the Royal Statistical Society. Series B
(Methodological) 32:432-434.\bigskip

Solari, M. (1969) The "Maximum Likelihood Solution" of the problem of
estimating a linear functional relationship. Journal of the Royal
Statistical Society. Series B (Methodological) 31:372-375.\bigskip

\end{document}